\documentclass[12pt]{article}
\usepackage{amsmath}
\usepackage{amsfonts}
\usepackage{geometry}
\usepackage[final]{graphicx}
\usepackage{float}
\usepackage[caption = false]{subfig}

\input{tcilatex}
\begin{document}

\title{\textbf{Certain Geometric Properties of Some Bessel Functions}}
\author{Rabha M. El-Ashwah$^{\text{1}}$\ and Alaa H. El-Qadeem$^{\text{2}}$
\and $^{\text{1}}${\small Department of Mathematics, Faculty of Science,
Damietta University, New Damietta 34517, Egypt} \and $^{\text{2}}${\small %
Department of Mathematics, Faculty of Science, Zagazig University, Zagazig
44519, Egypt} \and {\small r\_elashwah@yahoo.com} \\
{\small alaahassan1986@yahoo.com \& alaahassan@zu.edu.eg}}
\date{}
\maketitle

\begin{abstract}
In this paper, we determine necessary and sufficient conditions for the
generalized Bessel function to be in certain subclasses of starlike and
convex functions. Also, we obtain several corollaries as special cases of
the main results, these corollaries gives the corresponding results of the
familiar, modified and spherical Bessel functions.
\end{abstract}



\section*{1. Introduction}

Let $\mathcal{A}$ denotes the class of analytic functions which are defined
on the open unit disc $U=\{z\in
\mathbb{C}
:|z|$ $<1\}$, normalized by $f(0)=f^{\prime }(0)-1=0,$\ of the form%
\begin{equation}
f(z)=z+\dsum\limits_{k=2}^{\infty }a_{k}z^{k},  \tag{1.1}
\end{equation}%
\noindent For $0\leq \alpha <1$ and $0<\beta \leq 1,$\ let $S^{\ast }(\alpha
,\beta )$ denotes the subclass of $\mathcal{A}$ consisting of functions $%
f(z) $ of the form (1.1) and satisfy%
\begin{equation}
\left\vert \frac{\frac{zf^{\prime }(z)}{f(z)}-1}{\frac{zf^{\prime }(z)}{f(z)}%
+\left( 1-2\alpha \right) }\right\vert <\beta ;\text{ }z\in U,  \tag{1.2}
\end{equation}%
and $f\in K(\alpha ,\beta )$ denotes the subclass of $\mathcal{A}$
consisting of functions $f(z)$ of the form (1.1) such that%
\begin{equation}
\left\vert \frac{\frac{zf^{\prime \prime }(z)}{f^{\prime }(z)}}{\frac{%
zf^{\prime \prime }(z)}{f^{\prime }(z)}+2\left( 1-\alpha \right) }%
\right\vert <\beta ;\text{ }z\in U.  \tag{1.3}
\end{equation}

\noindent The subclasses $S^{\ast }(\alpha ,\beta )$ and $K(\alpha ,\beta )$
are the well-known subclasses of starlike and convex functions of order $%
\alpha $ and type $\beta $, respectively, introduced by Gupta and Jain \cite%
{Gupta}.\newline
\noindent It is clear that:%
\begin{equation}
f(z)\in K(\alpha ,\beta )\Longleftrightarrow zf^{\prime }(z)\in S^{\ast
}(\alpha ,\beta ).  \tag{1.4}
\end{equation}%
\noindent We note that $S^{\ast }(\alpha ,1)=S^{\ast }(\alpha )$ and $%
K(\alpha ,1)=K(\alpha )$ are the subclasses of starlike and convex functions
of order $\alpha (0\leq \alpha <1),$ respectively, also, $S^{\ast
}(0,1)=S^{\ast }(0)=S^{\ast }$ and $K(0,1)=K$ are the subclasses of starlike
and convex functions. These subclasses were introduced by Robertson \cite%
{Robertson}, Schild \cite{Schild}, MacGregor \cite{Mac}.

\bigskip

\noindent Recently, Baricz \cite{Baricz2008} defined a generalized Bessel
function $w_{p,b,c}(z)\equiv w(z)$ as follows:%
\begin{equation}
w_{p,b,c}(z)=\dsum\limits_{k=0}^{\infty }\frac{(-c)^{k}}{k!\Gamma \left( k+p+%
\frac{b+1}{2}\right) }\left( \frac{z}{2}\right) ^{2k+p}\left( p,b,c\in
\mathbb{C}
;\text{ }p+\frac{b+1}{2}\neq 0,-1,-2,...\right) ,  \tag{1.5}
\end{equation}%
which is the particular solution of the following second-order linear
homogeneous differential equation:%
\begin{equation}
z^{2}w^{\prime \prime }(z)+bzw^{\prime }(z)+\left[ cz^{2}-p^{2}+(1-b)p\right]
w(z)=0\text{ \ }\left( p,b,c\in
\mathbb{C}
\right) ,  \tag{1.6}
\end{equation}%
which, in turn, is a natural generalization of the classical Bessel's
equation. Solutions of (1.6) are regarded as the generalized Bessel function
of order $p$. The particular solution given by (1.5) is called the
generalized Bessel function of the first kind of order $p$.\newline

\noindent The series (1.5) permits the study of Bessel function, modified
Bessel functions and the spherical Bessel functions in a unified manner, as
follows:

\noindent \textit{(i)} Taking $b=c=1$ in (1.5), we have the familiar Bessel
function of the first kind of order $p$ defined by (see \cite{Watson} and
\cite{Baricz 2010})%
\begin{equation}
J_{p}(z)=\dsum\limits_{k=0}^{\infty }\frac{(-1)^{k}}{k!\Gamma \left(
p+k+1\right) }\left( \frac{z}{2}\right) ^{2k+p};\text{ }z\in
\mathbb{C}
.  \tag{1.7}
\end{equation}%
\noindent \textit{(ii)} Taking $b=1$ and $c=-1$ in (1.5), we obtain the
modified Bessel function of the first kind of order $p$ defined by (see \cite%
{Watson} and \cite{Baricz 2010})%
\begin{equation}
I_{p}(z)=\dsum\limits_{k=0}^{\infty }\frac{1}{k!\Gamma \left( p+k+1\right) }%
\left( \frac{z}{2}\right) ^{2k+p};\text{ }z\in
\mathbb{C}
.  \tag{1.8}
\end{equation}

\noindent \textit{(iii)} Taking $b=2$ and $c=1$ in (1.5), the function $%
w_{b,c,p}(z)$\ reduces to $\sqrt{2}j_{p}(z)/\sqrt{\pi },$\ where $j_{p}$\ is
the spherical Bessel function of the first kind of order $p$ defined by (see
\cite{Baricz 2010})%
\begin{equation}
j_{p}(z)=\sqrt{\frac{\pi }{2}}\dsum\limits_{k=0}^{\infty }\frac{(-1)^{k}}{%
k!\Gamma \left( p+k+\frac{3}{2}\right) }\left( \frac{z}{2}\right) ^{2k+p};%
\text{ }z\in
\mathbb{C}
.  \tag{1.9}
\end{equation}

\noindent Now, we will consider the function $u_{p,b,c}(z)$, which is the
normalized form of $w_{b,c,p}(z)$, as following:

\begin{equation}
u_{p,b,c}\left( z\right) =\Gamma \left( p+\tfrac{b+1}{2}\right) \text{ }z^{1-%
\frac{p}{2}}\text{ }w_{p,b,c}(2\sqrt{z}).  \tag{1.10}
\end{equation}%
By using the well-known Pochhammer symbol (or the shifted factorial)
defined, in terms of the familiar Gamma function, by%
\begin{equation*}
\left( \lambda \right) _{\mu }=\frac{\Gamma \left( \lambda +\mu \right) }{%
\Gamma \left( \lambda \right) }=\left\{
\begin{array}{c}
1,\text{ \ \ \ \ \ \ \ \ \ \ \ \ \ \ \ \ \ \ \ \ \ \ \ \ \ \ \ \ \ \ \ \ \ \
}\left( \mu =0\right) , \\
\lambda \left( \lambda +1\right) ...\left( \lambda +\mu -1\right) ,\text{ \
\ \ \ \ }\left( \mu \in
\mathbb{N}
\right) ,%
\end{array}%
\right.
\end{equation*}%
we can express $u_{p,b,c}(z)$ as follows:%
\begin{equation}
u_{p,b,c}(z)=\dsum\limits_{k=0}^{\infty }\frac{(-c)^{k}}{\left( p+\frac{b+1}{%
2}\right) _{k}}\frac{z^{k+1}}{k!},  \tag{1.11}
\end{equation}%
where $p+\frac{b+1}{2}\in
\mathbb{C}
\backslash
\mathbb{Z}
_{0}^{-}(%
\mathbb{Z}
_{0}^{-}:=%
\mathbb{Z}
^{-}\cup \{0\},%
\mathbb{Z}
^{-}=\left\{ -1,-2,...\right\} )$. Also, let $u_{p}^{\left( 1\right) }(z),$ $%
u_{p}^{\left( 2\right) }(z)$ and $u_{p}^{\left( 3\right) }(z)$ be the
normalized forms of $J_{p}(z)$, $I_{p}(z)$ and $j_{p}(z)$, respectively,
which are defined as following:%
\begin{equation}
u_{p}^{\left( 1\right) }(z):=\Gamma \left( p+1\right) z^{1-\frac{p}{2}}\text{
}J_{p}(2\sqrt{z})=\dsum\limits_{k=0}^{\infty }\frac{(-1)^{k}}{\left(
p+1\right) _{k}}\frac{z^{k+1}}{k!};\text{ }p\in
\mathbb{C}
\backslash
\mathbb{Z}
^{-},  \tag{1.11.a}
\end{equation}%
\begin{equation}
u_{p}^{\left( 2\right) }(z):=\Gamma \left( p+1\right) \text{ }z^{1-\frac{p}{2%
}}\text{ }I_{p}(2\sqrt{z})=\dsum\limits_{k=0}^{\infty }\frac{1}{\left(
p+1\right) _{k}}\frac{z^{k+1}}{k!};p\in
\mathbb{C}
\backslash
\mathbb{Z}
^{-},  \tag{1.11.b}
\end{equation}%
and%
\begin{equation}
u_{p}^{\left( 3\right) }(z):=\Gamma \left( p+\tfrac{3}{2}\right) \text{ }%
z^{1-\frac{p}{2}}\text{ }j_{p}(2\sqrt{z})=\dsum\limits_{k=0}^{\infty }\frac{%
(-1)^{k}}{\left( p+\frac{3}{2}\right) _{k}}\frac{z^{k+1}}{k!};p+\frac{1}{2}%
\in
\mathbb{C}
\backslash
\mathbb{Z}
^{-}.  \tag{1.11.c}
\end{equation}%
For further result on this transformation of the generalized Bessel
function, we refer to the recent papers \cite{ba2007}-\cite{bp} and \cite%
{dos}, where some interesting functional inequalities, integral
representations, extensions of some known trigonometric inequalities,
starlikeness, convexity and univalence, were established. Recently, Baricz
and Frasin \cite{bf} and Deniz et al. \cite{dos} were interested in the
univalence of some integral operators which involved the normalized form of
the ordinary Bessel function of the first kind and the normalized form of
the generalized Bessel functions of the first kind, respectively. Frasin
\cite{f2010} obtained various sufficient conditions for the convexity and
strong convexity of the integral operators defined by the normalized form of
the ordinary Bessel function of the first kind. Also, the problem of
geometric properties (such as univalence, starlikeness and convexity) of
some generalized integral operators has been discussed by many authors (see
\cite{bpp}, \cite{br}-\cite{bg}, \cite{f2008}, \cite{mb}, \cite{silverman2}).
For further properties we refer to El-Ashwah and Hassan \cite{H1}-\cite{H4}.

\noindent Very recently, Cho et al. \cite{clr} and Murugusundaramoorthy and
Janani \cite{mj} (see also Porwal and Dixit \cite{pd}) introduced some
characterization of generalized Bessel functions of first kind to be in
certain subclasses of uniformly starlike and uniformly convex functions.%
\newline
Motivated by the new technique due to Mustafa \cite{Nizami Mustafa},\ we
determine some geometric properties of generalized Bessel function of the
first type to be in $S^{\ast }(\alpha ,\beta )$ and $K(\alpha ,\beta )$,
followed by sufficient conditions for the familiar, modified and spherical
Bessel functions as special cases.

\section*{\protect\large 2. Main Result}

Unless otherwise mentioned, we assume in the reminder of this paper that$,$ $%
0\leq \alpha <1,$ $0<\beta \leq 1,$ $p,b,c\in
\mathbb{C}
,$ $p+\left( b+1\right) /2\in
\mathbb{C}
\backslash
\mathbb{Z}
_{0}^{-}$ and $z\in U$.\newline
To establish our main results, we shall require the following lemmas:

\noindent \textbf{Lemma 1} (\cite{Gupta}). \textit{A sufficient condition
for a function }$f$\textit{\ of the form (1.1) to be in the class }$S^{\ast
}(\alpha ,\beta )$ is that%
\begin{equation}
\dsum\limits_{k=2}^{\infty }\left[ k-1+\beta \left( k+1-2\alpha \right) %
\right] \left\vert a_{k}\right\vert \leq 2\beta \left( 1-\alpha \right) .
\tag{2.1}
\end{equation}

\noindent \textbf{Lemma 2} (\cite{Gupta}). \textit{A sufficient condition
for a function }$f$\textit{\ of the form (1.1) to be in the class }$K(\alpha
,\beta )$ is that%
\begin{equation}
\dsum\limits_{k=2}^{\infty }k\left[ k-1+\beta \left( k+1-2\alpha \right) %
\right] \left\vert a_{k}\right\vert \leq 2\beta \left( 1-\alpha \right) .
\tag{2.2}
\end{equation}%
\textbf{Theorem 1}. \textit{Let }$c<0$\textit{\ and }$p+\frac{b+1}{2}>0$%
\textit{\ and assume the following condition is satisfied}%
\begin{equation}
2\beta \left( 1-\alpha \right) \left[ 2-e^{^{\left( \frac{-c}{p+\frac{b+1}{2}%
+1}\right) }}+\frac{1}{p+\frac{b+1}{2}}\left( 1-e^{^{\left( \frac{-c}{p+%
\frac{b+1}{2}+1}\right) }}\right) \right] -\frac{\left( 1+\beta \right) (-c)%
}{\left( p+\frac{b+1}{2}\right) }e^{^{\left( \frac{-c}{p+\frac{b+1}{2}+1}%
\right) }}\geq 0,  \tag{2.3}
\end{equation}%
\textit{then }$u_{p,b,c}(z)\in S^{\ast }(\alpha ,\beta )$\textit{.}

\noindent \textit{Proof.} Since%
\begin{equation}
u_{p,b,c}(z)=z+\dsum\limits_{k=2}^{\infty }\frac{(-c)^{k-1}}{\left( p+\frac{%
b+1}{2}\right) _{k-1}}\frac{z^{k}}{(k-1)!},  \tag{2.4}
\end{equation}

\noindent by virtue of Lemma 1, it is suffices to show that
\begin{equation}
\dsum\limits_{k=2}^{\infty }\left[ k-1+\beta \left( k+1-2\alpha \right) %
\right] \frac{(-c)^{k-1}}{\left( p+\frac{b+1}{2}\right) _{k-1}(k-1)!}\leq
2\beta (1-\alpha ).  \tag{2.5}
\end{equation}%
\noindent Let%
\begin{equation*}
\tciLaplace (\alpha ,\beta ;p,b,c)=\dsum\limits_{k=2}^{\infty }\left[
k\left( 1+\beta \right) +\beta \left( 1-2\alpha \right) -1\right] \frac{%
(-c)^{k-1}}{\left( p+\frac{b+1}{2}\right) _{k-1}(k-1)!},
\end{equation*}%
by simple computation, we get%
\begin{eqnarray}
\tciLaplace (\alpha ,\beta ;p,b,c) &=&\dsum\limits_{k=2}^{\infty }\left[
\left( k-1\right) \left( 1+\beta \right) +2\beta \left( 1-\alpha \right) %
\right] \frac{(-c)^{k-1}}{\left( p+\frac{b+1}{2}\right) _{k-1}(k-1)!}  \notag
\\
&=&\left( 1{\small +}\beta \right) \dsum\limits_{k=2}^{\infty }\frac{%
(-c)^{k-1}}{\left( p+\frac{b+1}{2}\right) _{k-1}(k-2)!}+2\beta \left( 1%
{\small -}\alpha \right) \dsum\limits_{k=2}^{\infty }\frac{(-c)^{k-1}}{%
\left( p+\frac{b+1}{2}\right) _{k-1}(k-1)!}.  \TCItag{2.6}
\end{eqnarray}%
But%
\begin{eqnarray*}
\left( p+\tfrac{b+1}{2}\right) _{k-1} &=&\left( p+\tfrac{b+1}{2}\right)
\left( p+\tfrac{b+1}{2}+1\right) ...\left( p+\tfrac{b+1}{2}+k-2\right) \\
&\geq &\left( p+\tfrac{b+1}{2}\right) \left( p+\tfrac{b+1}{2}+1\right)
^{k-2};\text{ }k\geq 2,
\end{eqnarray*}%
which is equivalent to%
\begin{equation}
\frac{1}{\left( p+\frac{b+1}{2}\right) _{k-1}}\leq \left( p+\tfrac{b+1}{2}%
\right) \left( p+\tfrac{b+1}{2}+1\right) ^{k-2},  \tag{2.7}
\end{equation}%
then%
\begin{eqnarray*}
\tciLaplace (\alpha ,\beta ;p,b,c) &\leq &\frac{\left( 1+\beta \right) }{%
\left( p+\frac{b+1}{2}\right) }\dsum\limits_{k=2}^{\infty }\frac{(-c)^{k-1}}{%
\left( p+\frac{b+1}{2}+1\right) ^{k-2}(k-2)!} \\
&&+\frac{2\beta \left( 1-\alpha \right) }{\left( p+\frac{b+1}{2}\right) }%
\dsum\limits_{k=2}^{\infty }\frac{(-c)^{k-1}}{\left( p+\frac{b+1}{2}%
+1\right) ^{k-2}(k-1)!} \\
&=&\frac{\left( 1+\beta \right) (-c)}{\left( p+\frac{b+1}{2}\right) }%
\dsum\limits_{k=2}^{\infty }\frac{(-c)^{k-2}}{\left( p+\frac{b+1}{2}%
+1\right) ^{k-2}(k-2)!} \\
&&+\frac{2\beta \left( 1-\alpha \right) \left( p+\frac{b+1}{2}+1\right) }{%
\left( p+\frac{b+1}{2}\right) }\dsum\limits_{k=2}^{\infty }\frac{(-c)^{k-1}}{%
\left( p+\frac{b+1}{2}+1\right) ^{k-1}(k-1)!} \\
&=&\frac{\left( 1+\beta \right) (-c)}{\left( p+\frac{b+1}{2}\right) }%
e^{^{\left( \frac{-c}{p+\frac{b+1}{2}+1}\right) }}+\frac{2\beta \left(
1-\alpha \right) \left( p+\frac{b+1}{2}+1\right) }{\left( p+\frac{b+1}{2}%
\right) }\left( e^{^{\left( \frac{-c}{p+\frac{b+1}{2}+1}\right) }}-1\right) .
\end{eqnarray*}

\noindent It can be varified that the last expression $\tciLaplace (\alpha
,\beta ;p,b,c)$ is bounded above by $2\beta (1-\alpha )$ if (2.3) is
satisfied. This completes the proof of Theorem 1.

\bigskip

\noindent Taking $b=c=1$ in Theorem 1, we have the following result:

\noindent \textbf{Corollary 1}. \textit{If }$p>-1,$\textit{\ then the
condition}%
\begin{equation*}
2\beta \left( 1-\alpha \right) \left( e^{^{\frac{1}{p+2}}}\left( 2p+3\right)
-\left( p+2\right) \right) +\beta +1\geq 0,
\end{equation*}%
\textit{suffices that }$u_{p}^{\left( 1\right) }(z)\in S^{\ast }(\alpha
,\beta )$\textit{.}

\bigskip

\noindent Putting $b=1$ and $c=-1$\ in Theorem 1, we have the following
result:

\noindent \textbf{Corollary 2}. \textit{If }$p>-1,$\textit{\ then the
condition}%
\begin{equation*}
2\beta \left( 1-\alpha \right) \left( \left( 2p+3\right) -\left( p+2\right)
e^{^{\frac{1}{p+2}}}\right) +\left( \beta +1\right) e^{^{\frac{1}{p+2}}}\geq
0,
\end{equation*}%
\textit{\ suffices that }$u_{p}^{\left( 2\right) }(z)\in S^{\ast }(\alpha
,\beta )$\textit{.}\newline

\bigskip

\noindent Also, taking $b=2$ and $c=1$\ in Theorem 1, we have the following
result:

\noindent \textbf{Corollary 3}. \textit{If }$p>-\frac{3}{2},$\textit{\ then
the condition}%
\begin{equation*}
\beta \left( 1-\alpha \right) \left( 2p+4\right) \left( 2e^{^{\left( \frac{2%
}{2p+5}\right) }}-1\right) +\alpha \beta +1\geq 0,
\end{equation*}%
\textit{\ suffices that }$u_{p}^{\left( 3\right) }(z)\in S^{\ast }(\alpha
,\beta )$\textit{.}

\bigskip

\noindent Taking $\beta =1$ in Corollary 1, we have the following result:

\noindent \textbf{Corollary 4}. \textit{Let }$p>-1,$\textit{\ then }$%
u_{p}^{\left( 1\right) }(z)\in S^{\ast }(\alpha )$ if the following
condition is satisfied%
\begin{equation}
\left( 1-\alpha \right) \left( e^{^{\frac{1}{p+2}}}\left( 2p+3\right)
-\left( p+2\right) \right) +1\geq 0,  \tag{2.8}
\end{equation}

\noindent Taking $\alpha =0$ in Corollary 4, we have the following result:

\noindent \textbf{Corollary 5}. \textit{The function }$u_{p}^{\left(
1\right) }(z)$ is a starlike function for all $p>-1$.\newline
Proof. Taking $\alpha =0$ in (2.8) and using $p>-1,$ then the sufficient
condition for starlikeness of $u_{p}^{\left( 1\right) }(z)$ is that%
\begin{equation*}
\left( 2p+3\right) e^{^{\frac{1}{p+2}}}-\left( p+1\right) \geq 0,
\end{equation*}%
but the function $\left( 2x+3\right) e^{^{\frac{1}{x+2}}}-\left( x+1\right) $
it is non-negative for all $x>x_{0}\simeq -1.5314$ (see Figure 1). This
proves Corollary 5.
\begin{figure}[]
\centering
\includegraphics[width = 2in]{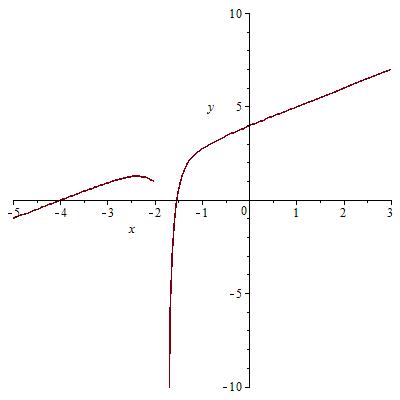}\newline
Figure 1
\end{figure}

\noindent \textbf{Remark 1.} The result introduced in Corollary 1 gives a
more improved result of this obtained by Mustafa \cite[Corollary 3.2]{Nizami
Mustafa} and the result obtained by Prajapat \cite[Corollary 2.8 (a)]%
{Prajapat2015}.

\noindent Taking $\beta =1$ in Corollary 2, we have the following result:

\noindent \textbf{Corollary 6}. \textit{Let }$p>-1,$\textit{\ then }$%
u_{p}^{\left( 2\right) }(z)\in S^{\ast }(\alpha )$ if the following
condition is satisfied%
\begin{equation}
\left( 1-\alpha \right) \left( \left( 2p+3\right) -\left( p+2\right) e^{^{%
\frac{1}{p+2}}}\right) +e^{^{\frac{1}{p+2}}}\geq 0,  \tag{2.9}
\end{equation}%
\noindent Taking $\alpha =0$ in Corollary 6, we have the following result:

\noindent \textbf{Corollary 7}. \textit{The function }$u_{p}^{\left(
2\right) }(z)$ is a starlike function for all $p>-1$.\newline
Proof. Taking $\alpha =0$ in (2.9) and using $p>-1,$ then the sufficient
condition for starlikeness of $u_{p}^{\left( 2\right) }(z)$ is that%
\begin{equation*}
\left( 2p+3\right) -e^{^{\frac{1}{p+2}}}\left( p+1\right) \geq 0,
\end{equation*}%
but the function $\left( 2x+3\right) -e^{^{\frac{1}{x+2}}}\left( x+1\right) $
it is non-negative for all $x>x_{0}=-2$ (see Figure 2). This proves
Corollary 7.
\begin{figure}[]
\centering
\includegraphics[width = 2in]{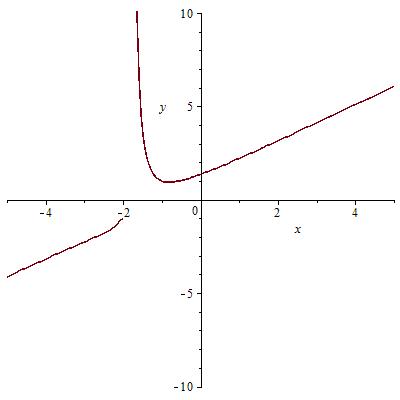}\newline
Figure 2
\end{figure}

\noindent Taking $\beta =1$ in Corollary 3, we have the following result:

\noindent \textbf{Corollary 8}. \textit{Let }$p>-\frac{3}{2},$\textit{\ then
}$u_{p}^{\left( 3\right) }(z)\in S^{\ast }(\alpha )$ if the following
condition is satisfied%
\begin{equation}
\left( 1-\alpha \right) \left( 2p+4\right) \left( 2e^{^{\left( \frac{2}{2p+5}%
\right) }}-1\right) +\alpha +1\geq 0,  \tag{2.10}
\end{equation}

\noindent Taking $\alpha =0$ in Corollary 8, we have the following result:

\noindent \textbf{Corollary 9}. \textit{The function }$u_{p}^{\left(
3\right) }(z)$ is a starlike function for all $p>-\frac{3}{2}$.\newline
Proof. Taking $\alpha =0$ in (2.10) and using $p>-\frac{3}{2},$ then the
sufficient condition for starlikeness of $u_{p}^{\left( 3\right) }(z)$ is
that%
\begin{equation*}
4\left( p+2\right) e^{^{\left( \frac{2}{2p+5}\right) }}-\left( 2p+3\right)
\geq 0,
\end{equation*}%
but the function $4\left( x+2\right) e^{^{\left( \frac{2}{2x+5}\right)
}}-\left( 2x+3\right) $ it is non-negative for all $x>x_{0}\simeq -2.0314$
(see Figure 3). This proves Corollary 9.
\begin{figure}[]
\centering
\includegraphics[width = 2in]{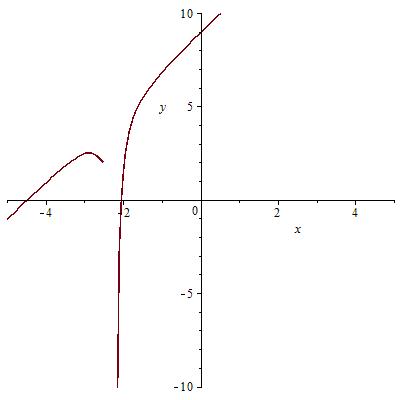}\newline
Figure 3
\end{figure}

\bigskip

\noindent \textbf{Theorem 2}. \textit{Let }$c<0$\textit{\ and }$p+\frac{b+1}{%
2}>0$\textit{\ and assume that the condition}%
\begin{equation}
2\beta (1{\small -}\alpha )\left( 1{\small +}\tfrac{p+\tfrac{b+1}{2}+1}{p+%
\tfrac{b+1}{2}}\right) e^{^{\frac{c}{p+\frac{b+1}{2}+1}}}{\small -}\left(
\tfrac{\left( 1+\beta \right) c^{2}}{(p+\tfrac{b+1}{2})(p+\tfrac{b+1}{2}+1)}%
{\small -}\tfrac{2\left( 1+\beta \left( 2-\alpha \right) \right) c}{p+\tfrac{%
b+1}{2}}{\small +}\tfrac{2\beta \left( 1-\alpha \right) (p+\tfrac{b+1}{2}+1)%
}{p+\tfrac{b+1}{2}}\right) {\small \geq }0,  \tag{2.11}
\end{equation}%
\textit{is satisfied, then }$u_{p,b,c}(z)\in K(\alpha ,\beta )$\textit{.}

\noindent \textit{Proof.} By virtue of Lemma 2, it is suffices to show that
\begin{equation}
\dsum\limits_{k=2}^{\infty }k\left[ k-1+\beta \left( k+1-2\alpha \right) %
\right] \frac{(-c)^{k-1}}{\left( p+\frac{b+1}{2}\right) _{k-1}(k-1)!}\leq
2\beta (1-\alpha ).  \tag{2.12}
\end{equation}%
\noindent Let%
\begin{eqnarray*}
\tciFourier (\alpha ,\beta ;p,b,c) &=&\dsum\limits_{k=2}^{\infty }k\left[
k\left( 1+\beta \right) +\beta \left( 1-2\alpha \right) -1\right] \frac{%
(-c)^{k-1}}{\left( p+\frac{b+1}{2}\right) _{k-1}(k-1)!} \\
&=&\dsum\limits_{k=2}^{\infty }\left[ k^{2}\left( 1+\beta \right) +k\left(
\beta \left( 1-2\alpha \right) -1\right) \right] \frac{(-c)^{k-1}}{\left( p+%
\frac{b+1}{2}\right) _{k-1}(k-1)!}
\end{eqnarray*}%
substituting $k^{2}=(k-1)(k-2)+3(k-1)+1$, $k=(k-1)+1$ and using (2.7), then
we have%
\begin{eqnarray*}
\tciFourier (\alpha ,\beta ;p,b,c) &=&\left( 1+\beta \right)
\dsum\limits_{k=2}^{\infty }\tfrac{(-c)^{k-1}}{\left( p+\frac{b+1}{2}\right)
_{k-1}(k-3)!}+\left( 2\beta \left( 2-\alpha \right) +2\underset{}{}\right)
\dsum\limits_{k=2}^{\infty }\tfrac{(-c)^{k-1}}{\left( p+\frac{b+1}{2}\right)
_{k-1}(k-2)!} \\
&&+\left( 2\beta \left( 1-\alpha \right) \underset{}{}\right)
\dsum\limits_{k=2}^{\infty }\tfrac{(-c)^{k-1}}{\left( p+\frac{b+1}{2}\right)
_{k-1}(k-1)!} \\
&\leq &\frac{\left( 1+\beta \right) (-c)^{2}}{\left( p+\tfrac{b+1}{2}\right)
\left( p+\tfrac{b+1}{2}+1\right) }\dsum\limits_{k=2}^{\infty }\tfrac{%
(-c)^{k-3}}{\left( p+\tfrac{b+1}{2}+1\right) ^{k-3}(k-3)!} \\
&&+\frac{\left( 2\beta \left( 2-\alpha \right) +2\right) (-c)}{\left( p+%
\tfrac{b+1}{2}\right) }\dsum\limits_{k=2}^{\infty }\tfrac{(-c)^{k-2}}{\left(
p+\tfrac{b+1}{2}+1\right) ^{k-2}(k-2)!} \\
&&+\frac{\left( 2\beta \left( 1-\alpha \right) \right) \left( p+\tfrac{b+1}{2%
}+1\right) }{\left( p+\tfrac{b+1}{2}\right) }\dsum\limits_{k=2}^{\infty }%
\tfrac{(-c)^{k-1}}{\left( p+\tfrac{b+1}{2}+1\right) ^{k-1}(k-1)!} \\
&=&\frac{\left( 1+\beta \right) (-c)^{2}}{\left( p+\tfrac{b+1}{2}\right)
\left( p+\tfrac{b+1}{2}+1\right) }e^{^{\left( \frac{-c}{p+\frac{b+1}{2}+1}%
\right) }}+\frac{\left( 2\beta \left( 2-\alpha \right) +2\right) (-c)}{%
\left( p+\tfrac{b+1}{2}\right) }e^{^{\left( \frac{-c}{p+\frac{b+1}{2}+1}%
\right) }} \\
&&+\frac{\left( 2\beta \left( 1-\alpha \right) \right) \left( p+\tfrac{b+1}{2%
}+1\right) }{\left( p+\tfrac{b+1}{2}\right) }\left( e^{^{\left( \frac{-c}{p+%
\frac{b+1}{2}+1}\right) }}-1\right) \\
&=&\left[ \tfrac{\left( 1+\beta \right) (-c)^{2}}{\left( p+\tfrac{b+1}{2}%
\right) \left( p+\tfrac{b+1}{2}+1\right) }+\tfrac{\left( 2\beta \left(
2-\alpha \right) +2\right) (-c)}{\left( p+\tfrac{b+1}{2}\right) }+\tfrac{%
\left( 2\beta \left( 1-\alpha \right) \right) \left( p+\tfrac{b+1}{2}%
+1\right) }{\left( p+\tfrac{b+1}{2}\right) }\right] e^{^{\left( \frac{-c}{p+%
\frac{b+1}{2}+1}\right) }} \\
&&-\tfrac{\left( 2\beta \left( 1-\alpha \right) \right) \left( p+\tfrac{b+1}{%
2}+1\right) }{\left( p+\tfrac{b+1}{2}\right) }.
\end{eqnarray*}

\noindent The last expression $\tciFourier (\alpha ,\beta ;p,b,c)$ is
bounded above by $2\beta (1-\alpha )$ if (2.11) is satisfied. This completes
the proof of Theorem 2.

\bigskip

\noindent Taking $b=c=1$ in Theorem 2, we have the following result:

\noindent \textbf{Corollary 10}. \textit{If }$p>-1,$\textit{\ then the
condition}%
\begin{equation*}
2\beta (1{\small -}\alpha )\left( 1{\small +}\tfrac{p+2}{p+1}\right) e^{^{%
\frac{1}{p+2}}}{\small -}\left( \tfrac{1+\beta }{(p+1)(p+2)}{\small +}\tfrac{%
2\beta \left( 1-\alpha \right) (p+2)}{p+1}{\small -}\tfrac{2\left( 1+\beta
\left( 2-\alpha \right) \right) }{p+1}\right) {\small \geq }0,
\end{equation*}%
\textit{suffices that }$u_{p}^{\left( 1\right) }(z)\in K(\alpha ,\beta )$%
\textit{.}

\bigskip

\noindent Putting $b=1$ and $c=-1$\ in Theorem 2, we have the following
result:

\noindent \textbf{Corollary 11}. \textit{If }$p>-1,$\textit{\ then the
condition}%
\begin{equation*}
2\beta (1{\small -}\alpha )\left( 1{\small +}\tfrac{p+2}{p+1}\right) {\small %
-}\left( \tfrac{\left( 1+\beta \right) }{(p+1)(p+2)}{\small +}\tfrac{2\left(
1+\beta \left( 2-\alpha \right) \right) }{p+1}{\small +}\tfrac{2\beta \left(
1-\alpha \right) (p+2)}{p+1}\right) e^{^{\frac{1}{p+2}}}{\small \geq }0,
\end{equation*}%
\textit{\ suffices that }$u_{p}^{\left( 2\right) }(z)\in K(\alpha ,\beta )$%
\textit{.}\newline

\bigskip

\noindent Also, taking $b=2$ and $c=1$\ in Theorem 2, we have the following
result:

\noindent \textbf{Corollary 12}. \textit{If }$p>-\frac{3}{2},$\textit{\ then
the condition}%
\begin{equation*}
\beta (1{\small -}\alpha )\left( 1{\small +}\tfrac{2p+5}{2p+3}\right) e^{^{%
\frac{2}{2p+5}}}{\small -}\left( \tfrac{2\left( 1+\beta \right) }{%
(2p+3)(2p+5)}{\small +}\tfrac{\beta \left( 1-\alpha \right) (2p+5)}{2p+3}%
{\small -}\tfrac{2\left( 1+\beta \left( 2-\alpha \right) \right) }{2p+3}%
\right) {\small \geq }0,
\end{equation*}%
\textit{suffices that }$u_{p}^{\left( 3\right) }(z)\in K(\alpha ,\beta )$%
\textit{.}

\bigskip

\noindent Taking $\beta =1$ in Corollary 10, we have the following result:

\noindent \textbf{Corollary 13}. \textit{Let }$p>-1,$\textit{\ then }$%
u_{p}^{\left( 1\right) }(z)\in K(\alpha )$ if the following condition is
satisfied%
\begin{equation}
(1{\small -}\alpha )\left( 1{\small +}\tfrac{p+2}{p+1}\right) e^{^{\frac{1}{%
p+2}}}{\small -}\left( \tfrac{1}{(p+1)(p+2)}{\small +}\tfrac{\left( 1-\alpha
\right) (p+2)}{p+1}{\small -}\tfrac{3-\alpha }{p+1}\right) {\small \geq }0,
\tag{2.13}
\end{equation}%
\noindent Taking $\alpha =0$ in Corollary 13, we have the following result:

\noindent \textbf{Corollary 14}. \textit{The function }$u_{p}^{\left(
1\right) }(z)$ is a convex function for all $p>-1$.\newline
Proof. Taking $\alpha =0$ in (2.13) and using $p>-1,$ then the sufficient
condition for starlikeness of $u_{p}^{\left( 1\right) }(z)$ is that%
\begin{equation}
\left( 2p^{2}+7p+6\right) e^{^{\frac{1}{p+2}}}{\small -}\left(
p^{2}+p-1\right) {\small \geq }0,  \tag{2.13}
\end{equation}%
but the function $\left( 2x^{2}+7x+6\right) e^{^{\frac{1}{x+2}}}{\small -}%
\left( x^{2}+x-1\right) $ it is non-negative for all $x>x_{0}\simeq -1.5254$
(see Figure 4). This proves Corollary 14.
\begin{figure}[]
\centering
\includegraphics[width = 2in]{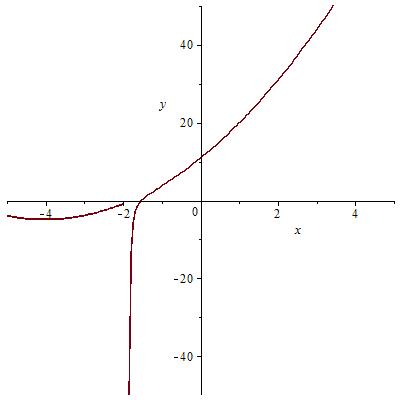}\newline
Figure 4
\end{figure}

\noindent \textbf{Remark 2.} The result introduced in Corollary 14 gives a
more improved result of this obtained by Mustafa \cite[Corollary 3.4]{Nizami
Mustafa}.

\noindent Taking $\beta =1$ in Corollary 11, we have the following result:

\noindent \textbf{Corollary 15}. \textit{Let }$p>-1,$\textit{\ then }$%
u_{p}^{\left( 2\right) }(z)\in K(\alpha )$ if the following condition is
satisfied%
\begin{equation}
(1{\small -}\alpha )\left( 1{\small +}\tfrac{p+2}{p+1}\right) {\small -}%
\left( \tfrac{1}{(p+1)(p+2)}{\small +}\tfrac{3-\alpha }{p+1}{\small +}\tfrac{%
\left( 1-\alpha \right) (p+2)}{p+1}\right) e^{^{\frac{1}{p+2}}}{\small \geq }%
0,  \tag{2.14}
\end{equation}

\noindent Taking $\alpha =0$ in Corollary 15, we have the following result:

\noindent \textbf{Corollary 16}. \textit{The function }$u_{p}^{\left(
2\right) }(z)$ is a convex function for all $p>x_{0}\simeq 3.8523$.\newline
Proof. Taking $\alpha =0$ in (2.14) and using $p>-1,$ then the sufficient
condition for starlikeness of $u_{p}^{\left( 2\right) }(z)$ is that%
\begin{equation*}
\left( 2p^{2}+7p+6\right) {\small -}\left( p^{2}+7p+11\right) e^{^{\frac{1}{%
p+2}}}{\small \geq }0,
\end{equation*}%
but the function $\left( 2x^{2}+7x+6\right) {\small -}\left(
x^{2}+7x+11\right) e^{^{\frac{1}{x+2}}}$ it is non-negative for all $%
x>x_{0}\simeq 3.8523$ (see Figure 5). This proves Corollary 16.

\begin{figure}[]
\centering
\includegraphics[width = 2in]{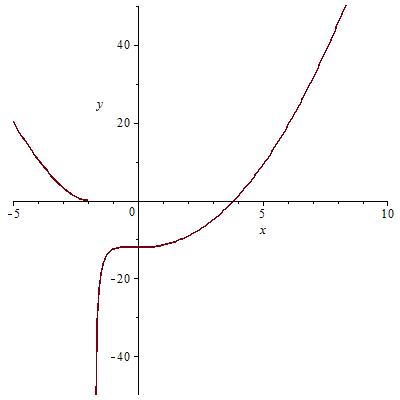}\newline
Figure 5
\end{figure}

\noindent Taking $\beta =1$ in Corollary 12, we have the following result:

\noindent \textbf{Corollary 17}. \textit{Let }$p>-\frac{3}{2},$\textit{\
then }$u_{p}^{\left( 3\right) }(z)\in K(\alpha )$ if the following condition
is satisfied%
\begin{equation}
(1{\small -}\alpha )\left( 1{\small +}\tfrac{2p+5}{2p+3}\right) e^{^{\frac{2%
}{2p+5}}}{\small -}\left( \tfrac{4}{(2p+3)(2p+5)}{\small +}\tfrac{\left(
1-\alpha \right) (2p+5)}{2p+3}{\small -}\tfrac{2\left( 3-\alpha \right) }{%
2p+3}\right) {\small \geq }0,  \tag{2.15}
\end{equation}

\noindent Taking $\alpha =0$ in Corollary 17, we have the following result:

\noindent \textbf{Corollary 18}. \textit{The function }$u_{p}^{\left(
3\right) }(z)$ is a convex function for all $p>-\frac{3}{2}$.\newline
Proof. Taking $\alpha =0$ in (2.15) and using $p>-\frac{3}{2},$ then the
sufficient condition for starlikeness of $u_{p}^{\left( 3\right) }(z)$ is
that%
\begin{equation*}
\left( 8p^{2}+36p+40\right) e^{^{\frac{2}{2p+5}}}{\small -}\left( 4p^{2}%
{\small +8p-1}\right) {\small \geq }0,
\end{equation*}%
but the function $\left( 8x^{2}+36x+40\right) e^{^{\frac{2}{2x+5}}}{\small -}%
\left( 4x^{2}+8x-1\right) $ it is non-negative for all $x>x_{0}\simeq
-2.0254 $ (see Figure 6). This proves Corollary 18.
\begin{figure}[]
\centering
\includegraphics[width = 2in]{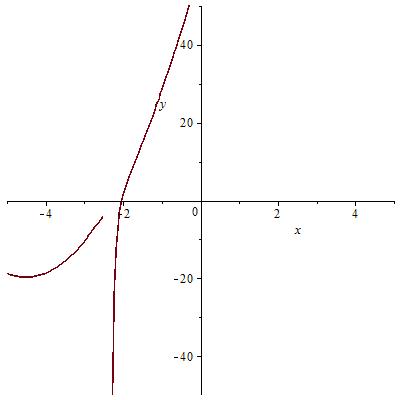}\newline
Figure 6
\end{figure}

\end{document}